\documentclass[oneside,english]{amsart}
\usepackage[T1]{fontenc}
\usepackage[latin9]{inputenc}
\usepackage{units}
\usepackage{amsthm}
\usepackage{amssymb}
\usepackage{esint}

\makeatletter
\numberwithin{equation}{section}
\numberwithin{figure}{section}
  \theoremstyle{definition}
  \newtheorem*{example*}{Example}
  \theoremstyle{remark}
  \newtheorem*{rem*}{Remark}
\theoremstyle{plain}
\newtheorem{thm}{Theorem}
  \theoremstyle{plain}
  \newtheorem{cor}[thm]{Corollary}
  \theoremstyle{definition}
  \newtheorem{defn}[thm]{Definition}
  \theoremstyle{plain}
  \newtheorem{prop}[thm]{Proposition}

\makeatother

\usepackage{babel}

\begin{document}

\title{Stability of the global attractor under Markov-Wasserstein noise}

\author{Martin Kell}

\date{March 17, 2011}

\address{\noindent \tiny Max-Planck-Institute for Mathematics in the Sciences,
Inselstr. 22-26, D-04103 Leipzig, Germany}

\email{mkell@mis.mpg.de}

\thanks{The author would like to thank the IMPRS {}``Mathematics in the
Sciences'' for financial support and his advisor, Prof. J\"urgen
Jost, and the MPI MiS for providing an inspiring research atmosphere.}

\subjclass[2000]{34D23, 37B35, 60B05, 60B10}

\keywords{global attractor, random perturbation, Wasserstein space}
\begin{abstract}
We develop a {}``weak Wa{\.z}ewski principle'' for discrete and
continuous time dynamical systems on metric spaces having a weaker
topology to show that attractors can be continued in a weak sense.
After showing that the Wasserstein space of a proper metric space
is weakly proper we give a sufficient and necessary condition such
that a continuous map (or semiflow) induces a continuous map (or semiflow)
on the Wasserstein space. In particular, if these conditions hold
then the global attractor, viewed as invariant measures, can be continued
under Markov-type random perturbations which are sufficiently small
w.r.t. the Wasserstein distance, e.g. any small bounded Markov-type
noise and Gaussian noise with small variance will satisfy the assumption.
\end{abstract}
\maketitle
In this paper we are going to show that the invariant measures of
a dynamical system having a global attractor (either discrete and
continuous time) can be {}``continued'' under small (not necessarily
bounded) noise. Instead of just showing that there is a stationary
measure of the perturbed system {}``weakly'' close to the original
ones (see e.g. {\cite[1.7]{Kifer1988}}) we show that it is close
w.r.t. the Wasserstein metric $w_{p}$ (the order $p$ depends on
regularity of the noise). Previous research mainly focused on Gaussian
type noise, {}``absolute continuous'' noise or assumed implicitly
bounded noise, e.g. Kifer {\cite[p.103]{Kifer1988}} and L.-S. Young
\cite{Young1986} considered noise on a positive invariant (bounded)
neighborhood $U$ of a local attractor which is zero on $\partial U$
and thus is bounded. If $X$ is compact than any noise will be bounded.
Thus we are in particular interested in non-compact $X$, although
we restrict our attention to proper metric spaces which includes all
locally compact geodesic spaces.

All our results apply equally to discrete and continuous time dynamical
system. We will mainly focus on discrete time because there is a better
intuition behind these. In the first section we extend Rybakowski's
continuation of a positive invariant isolating neighborhood to the
discrete time setting which will be the key step to treat continuous
and discrete time systems on equal footing. Using a kind of {}``weak
Wa{\.z}ewski principle'' we show that attractors can be continued
in a weak sense without assuming admissibility of the perturbed system
(theorem \ref{thm:weak-continuation}).

Then we introduce the Wasserstein space and show that the Wasserstein
space of a proper metric space is weakly proper, i.e. closed $\delta$-neighborhoods
of compact sets are weakly compact (theorem \ref{thm:weakly-proper}).
We give a necessary and sufficient condition such that a dynamical
system $f$ (resp. a semiflow $\pi$) on a proper metric space makes
the transfer map $f_{*}:\mathcal{P}(X)\to\mathcal{P}(X)$ (resp. transfer
semiflow $\pi_{*}$), which is always (weakly) continuous, (strongly)
continuous on any Wasserstein space $(\mathcal{P}_{p}(X),w_{p})$
of order $1\le p<\infty$. 

Finally, we look at Markov-type perturbations of a dynamical system
$f$ (resp. semiflow $\pi$) having a global attractor. These are
perturbation of $f_{*}$ in the Wasserstein space. Under quite general
assumptions on $f$ we can show using Conley theory that if the perturbation
is sufficiently small then the perturbed system has an isolated (weak)
attractor and there is at least one stationary measure (theorem \ref{thm:weak-cont-wasserstein}).
The perturbed attractor is strongly close to the unperturbed attractor
and contains all invariant measures and its positive invariant isolating
neighborhood is convex. If the noise is of order $p$ then the mass
of the invariant measures decays at least as $R^{-p}$ where $R$
is the distance from the global attractor of $f$ (resp. $\pi$).
Furthermore, a standard result from Conley theory shows that if the
noise level converges to zero then the set of invariant measures converges
in the Wasserstein space to the set of invariant measures of the deterministic
system (as usual without further assumptions only upper semicontinuity
holds).

Because of the special form of the Kantorovich-Rubinstein duality
in $\mathcal{P}_{1}(X)$ we can show that the (local) weak attractor
of the perturbed system $P:\mathcal{P}_{1}(X)\to\mathcal{P}_{1}(X)$
is actually the global weak attractor of $P$.

The framework of a metric space with a weaker topology used here is
similar to the framework in {\cite[section 2.1]{AmbGigSav2008}} used
to construct general gradient flows.

\subsection*{Motivation}

Consider a dynamical system $f:X\to X$ on a proper metric space $X$
having a global attractor, i.e. a compact invariant set $A$ that
attracts all its (bounded) neighborhoods. We will not consider the
map $f$ itself, but the map $f_{*}:\mathcal{P}_{1}(X)\to\mathcal{P}_{1}(X)$
defined via push-forward map on the space of probability measures.
If $f$ is {}``nice'' then $f_{*}$ is continuous and has the global
attractor $K=\mathcal{P}(A)$, i.e. the probability measures supported
on the global attractor of $f$. 

Markov-type noise can be considered as a perturbation $\tilde{F}$
of $f_{*}$, i.e. instead of $\delta_{x}\mapsto\delta_{f(x)}$ we
have $\delta_{x}\mapsto p(dy|x)$ where $p(dy|x)$ and $\delta_{f(x)}$
uniformly close w.r.t. the Wasserstein distance $w_{1}$ for all $x$.
This can be seen as a smearing of the image $f(x)$ or some uncertainty
about the actual image. For example, if $f$ is the time-$1$ map
of a flow generated by the ODE $\dot{x}=g(x)$ then $\tilde{F}$ could
be the time-one map of the flow of distributions of the SDE $dx=g(x)dt+\epsilon dW_{t}$,
i.e. additive Gaussian noise with small variance. 

If $\tilde{F}$ and $f_{*}$ are sufficiently close then $\tilde{F}$
has a (weak) global attractor (in $\mathcal{P}_{1}(X)$) which is
close to $K$ w.r.t. $w_{1}$. Hence stability of the global attractor
holds in the Wasserstein space $\mathcal{P}_{1}(X)$. 

The following example is inspired by Crauel, Flandoli - {}``Additive
Noise Destroys a Pitchfork Bifurcation'' \cite{Crauel1998} and could
be decribed as {}``Additive Noise Destroys Attractors''. The noise
will be worse than white noise used by Crauel and Flandoli, but can
still be considered as small. 
\begin{example*}
[Generic collapse under "small" noise](1) Suppose $f:X\to X$ has
a global attractor and at least one fixed point $x_{0}$ (the argument
works equally well with general attractors). Take any noise level
$\epsilon>0$ and let $P_{\epsilon}:\mathcal{P}(X)\to\mathcal{P}(X)$
be the Markov map induced by \[
x\mapsto(1-\epsilon)\delta_{f(x)}+\epsilon\delta_{x_{0}}.\]
This map is (weakly) close to the unperturbed system $f_{*}:\mathcal{P}(X)\to\mathcal{P}(X)$.
Namely, if $d_{LP}$ is the Levy-Prokhorov distance (which metrizes
$\mathcal{P}(X)$) then \[
\sup_{\mu\in P(X)}d_{LP}(P_{\epsilon}(\mu),f_{*}(\mu))\le\epsilon.\]
 But $P$ has exactly one invariant measure, namely $\delta_{x_{0}}$,
and all others converge to this measure weakly.

(2) Now we want to show that this can also happen in any Wasserstein
space $\mathcal{P}_{p}(X)$ for $1\le p<\infty$ ($\mathcal{P}_{\infty}(X)$
only allows bounded noise which when sufficiently small cannot destroy
local attractors and thus a global attractors with at least two sinks
never collapses, see \cite{Kell2011b}). Suppose $a:X\to[0,1]$ is
a continuous function. Define $q_{a}(dy|\cdot):X\to\mathcal{P}_{p}(X)$
by \[
x\mapsto q_{a}(dy|x)=(1-a(x))\delta_{f(x)}+a(x)\delta_{x_{0}},\]
which is obviously continuous. Thus \[
w_{p}(q_{a}(dy|x),\delta_{f(x)})^{p}=a(x)d(f(x),x_{0})^{p}.\]
 So if we define \[
a(x)=\frac{\epsilon^{p}}{1+d(f(x),x_{0})^{p}}\]
 then \[
w_{p}(q_{a}(dy|x),\delta_{f(x)})^{p}=\frac{\epsilon^{p}\cdot d(f(x),x_{0})^{p}}{1+d(f(x),x_{0})^{p}}\le\epsilon^{p}.\]
So, in particular, the induced MW-map $Q_{a}:\mathcal{P}_{p}(X)\to\mathcal{P}_{p}(X)$
of order $p$ relative to $f$ has noise level $\epsilon$. Furthermore,
the only invariant measure of $Q_{a}$ is $\delta_{x_{0}}$ and all
other measures converge to it.
\end{example*}
The example above should make clear that using arbitrary unbounded
noise even when it is small can have strange effects on the global
attractor. Although we have some {}``attracting'' invariant measures
of the perturbed system the attractor might look very different from
the original one, in our case it might be just one fixed point and
this one can even be the {}``most'' unstable one of the original
attractor. Therefore, stochastic stability of attractors under arbitrary
{}``small'' noise should not be referred to a single invariant measure
but to all of them, even though we can speak of stochastic stability
if the type of noise is more restricted, besides of being sufficiently
{}``small''.

\section{Discrete-time Conley theory for stable invariant sets}

In this section we will use Conley theory, that is continuation methods
from Conley index theory without using the topological (or (co)homological)
Conley index. We will prove a continuation for a positive invariant
neighborhood of a stable isolated invariant set of a time discrete
dynamical system. The result will not require a compactness assumption
(called admissibility) of the perturbed system and is a different
type of continuation than \cite{MroRyb1991}. Our proof will follow
the proof of {\cite[Theorem 12.3]{Rybakowski1987}} which is the continuation
for semiflows. In particular, the results stated here and in the next
sections also hold for semiflows if we assume that they do not explode
on a given neighborhood. 

In both cases the Wa{\.z}ewski principle for the index pair of the
perturbed system does not apply. But we can use other assumptions
to show that attractors continue, e.g. the map is weakly continuous
and closed $\delta$-neighborhoods of compact set are weakly compact,
which is the case for the Wasserstein space on proper metric spaces.

We will now give the definitions used in \cite{MroRyb1991} and \cite{Rybakowski1987}
to prove the existence of an index pair for certain isolated invariant
sets. Our setting will be a complete separable metric space $Y$ and
a dynamical system, i.e. a continuous map $f:Y\to Y$. A full left
solution of $f$ in $N$ is a sequence $\{x_{-n}\}_{n\in\mathbb{N}}\subset N$
such that $f(x_{n-1})=x_{n}$ for $n\le0$. Define the following sets\begin{eqnarray*}
A^{+}(N) & = & \{x\in N\,|\, f^{k}(x)\in N\,\mbox{for all }k\ge0\}\\
A^{-}(N) & = & \{x\in N\,|\,\exists\mbox{\,\ full left solution \ensuremath{\{x_{-n}\}_{n\in\mathbb{Z}}}\,\ in \,\ensuremath{N}\,\ through \ensuremath{x_{0}=x}}\}\\
A(N) & = & A^{+}(N)\cap A^{-}(N).\end{eqnarray*}
These are called the maximal positive invariant (resp. negative invariant,
resp. invariant) set in $N$. If $N$ is unbounded then $A(N)$ usually
denotes only the bounded invariant  orbits instead of all of them.
A set $K$ is called invariant if $A(K)=K$. If there is a closed
neighborhood $N$ of an invariant set $K$ with $A(N)=K$ then $K$
is called isolated with isolating neighborhood $N$.

For $l,m\in\mathbb{N}$ and $l\le k$ define \[
f^{[l,m]}(x)=\{y\,|\, f^{k}(x)=y\,\mbox{for some }k\in\mathbb{N}\cap[l,m]\}\]
If $f_{n}:Y\to Y$ is a sequence of continuous maps such that $f_{n}\to f$,
i.e. $f_{n}(x_{n})\to f(x)$ whenever $x_{n}\to x$ as $n\to\infty$,
then we say that a closed bounded set $N$ is $\{f_{n}\}$-admissible
if for any sequence $\{x_{n}\}_{n\in\mathbb{N}}$ with $f_{n}^{[0,m_{n}]}(x_{n})\subset N$
and $m_{n}\to\infty$ the sequence of endpoints $\{f_{n}^{m_{n}}(x_{n})\}_{n\in\mathbb{N}}$
is precompact. In case this property holds for $f_{n}\equiv f$ then
we just say $N$ is $f$-admissible.
\begin{rem*}
Later on, we deal with dynamical systems on the space of probability
measures $Y=\mathcal{P}(X)$ for some metric space $X$. An invariant
measure for that system is invariant w.r.t. the definition above.
In particular, periodic measures will be called invariant. A fixed
point for these systems will be called stationary measure. 
\end{rem*}
In the following we will use several ideas from \cite{MroRyb1991}:
Let $N,N'$ be two $f$-admissible isolating neighborhood for some
isolated invariant set $K$ with \[
N\subset\operatorname{int}N'\cap f^{-1}(\operatorname{int}N').\]
 The authors in {\cite[4.4]{MroRyb1991}} used a so called Lyapunov
pair $(\phi,\gamma)$ which is continuous on a small neighborhood
$W\subset N$ of $K$ and has the following properties: $K\subset\gamma^{-1}(0)$,
$\phi$ (resp. $\gamma$) is decreasing (resp. increasing) along orbits
and $\phi(x)=0$ with $x\in W$ implies \[
x\in A^{-}(N)\cup\partial N'.\]
 Because $K$ is compact we can choose $d(W,\partial N')>0$ and assume
$x\in A^{-}(N)$ whenever $\phi(x)=0$. Furthermore, it is shown that
if $\phi(x_{n})\to0$ then $x_{n}$ admits a convergent subsequence.
\begin{thm}
\label{thm:stable}Suppose $N$ is an isolating neighborhood for $K$
such that the assumptions above hold and \[
A^{-}(N)=A(N)=K\ne\varnothing.\]
Then there exists an admissible isolating neighborhood $B\subset N$
which is positive invariant, i.e. no trajectories exit $B$.\end{thm}
\begin{rem*}
This is the discrete time version of {\cite[I-5.5]{Rybakowski1987}}
using the theory of \cite{MroRyb1991}. The proof is essentially copied
from Rybakowski using the Lyapunov pair above.\end{rem*}
\begin{proof}
Define \begin{eqnarray*}
P_{1}^{\epsilon} & = & N\cap\operatorname{cl}\{x\in\operatorname{int}N'\,|\,\phi(x)<\epsilon\}\\
P_{2}^{\epsilon} & = & P_{1}^{\epsilon}\backslash\{x\in\operatorname{int}N'\,|\,\gamma(x)<\epsilon\}.\end{eqnarray*}
 it was shown {\cite[4.4]{MroRyb1991}} that $P_{1}^{\epsilon}\subset W$
is a neighborhood of $K$ for sufficiently small $\epsilon>0$ and
whenever $x\in P_{i}^{\epsilon}$ and $f(x)\in N$ then $x\in P_{i}^{\epsilon}$
and if $x\in P_{1}^{\epsilon}$ and $f(x)\notin N$ then $x\in P_{2}^{\epsilon}$,
i.e. $P_{2}^{\epsilon}$ is the exit ramp for $P_{1}^{\epsilon}$.

Now fix a sufficiently small $\epsilon>0$ and let $0<\delta\le\epsilon$
then $P_{1}^{\delta}\subset P_{1}^{\epsilon}$. Define $\tilde{P}_{2}^{\delta}:=P_{1}^{\delta}\cap P_{2}^{\epsilon}$
then because $\phi$ is decreasing along orbits $\tilde{P}_{2}^{\delta}$
is still an exit ramp for $P_{2}^{\delta}$.

If $A^{-}(N)=A(N)=K$ then we claim that there is a $\delta>0$ such
that $\tilde{P}_{2}^{\delta}=\varnothing$ which implies that $P_{1}^{\delta}$
is positive invariant. If this does not hold then there is a sequence
$x_{n}\in P_{1}^{\delta_{n}}\cap P_{2}^{\epsilon}$ with $\delta_{n}\to0$.
Thus $\phi(x_{n})\to0$ and $\gamma(x_{n})\ge\epsilon$ which implies
that there is a subsequence $x_{n'}\to x\in N$ such that $\phi(x)=0$.
Hence \[
x\in A^{-}(N)=K\]
and $\gamma(x)=0$ by assumption. But $\gamma$ is continuous and
$x_{n'}\to x$ implies $\epsilon\le\gamma(x_{n'})\to\gamma(x)=0$
which is a contradiction. This proofs our claim and thus the theorem.
\end{proof}
In the following we will assume that $K$ satisfies the assumption
of the theorem and that $B:=P_{1}^{\delta_{0}}$, $\phi$ and $\gamma$
are given as in the proof. It is obvious that $P_{1}^{\delta}$ is
positive invariant w.r.t. $f$ for any $0<\delta\le\delta_{0}$. Furthermore,
suppose $f_{n}\to f$. 
\begin{thm}
\label{thm:cont}Assume $N'$ (see above) is $\{f_{n_{m}}\}$-admissible
for each subsequence of $\{f_{n}\}_{n\in\mathbb{N}}$. Set $\tilde{U}=\operatorname{int}B$
and define \[
V(a)=\{x\in\tilde{U}\,|\,\phi(x)<a\}.\]
Then for some $a_{0}>0$, $N:=\operatorname{cl}V(a_{0})\subset\tilde{U}$.
Furthermore, for some sufficiently small $\epsilon_{0}>0$ and all
$0<\epsilon\le\epsilon_{0}$ there is an $n_{0}=n_{0}(\epsilon)$
such that for all $n\ge n_{0}$ there is a positive $f_{n}$-invariant
closed $N_{n}(\epsilon)$ and \[
K_{n}\subset V(\epsilon)\subset N_{n}(\epsilon)\subset N.\]
\end{thm}
\begin{rem*}
The complete continuation theorem for index pairs does not hold for
discrete time dynamical systems in general. A proof would require
that there is a neighborhood such that the exit time is continuous
in $\tilde{U}$ , i.e. $\omega_{n}^{+}(x_{n})\to\omega^{+}(x_{0})$
whenever $x_{n}\to x_{0}$ in $\tilde{U}$, which holds for semiflows
only for so called isolating blocks. These blocks do not necessarily
exist for continuous maps.\end{rem*}
\begin{proof}
By {\cite[3.9]{MroRyb1991}} there is an $a_{0}>0$ such that $N=\operatorname{cl}V(a_{0})\subset\tilde{U}$.
And similar to {\cite[I-4.5]{Rybakowski1987}} we can show that for
$0<\epsilon\le a_{0}$ and all $n\ge n_{0}(\epsilon)$ \[
K_{n}\subset V(\epsilon).\]

Define \[
N_{n}(\epsilon)=N\cap\operatorname{cl}\{y\,|\,\mbox{ \ensuremath{\exists x\in V(\epsilon)}, \ensuremath{m\ge0}\,\ s.t.\,\ \ensuremath{f_{n}^{[0,m]}(x)\subset\tilde{U}\,}and \ensuremath{f_{n}^{m}(x)=y}}\}.\]
 Following the proof of {\cite[I-12.5]{Rybakowski1987}} we can show
that $N_{n}(\epsilon)$ satisfies the following properties for $n\ge n_{0}(\epsilon)$
\begin{itemize}
\item $x\in N_{n}(\epsilon)$ and $f_{n}(x)\in N$ implies $f_{n}(x)\in N_{n}(\epsilon)$
\item $K_{n}\subset V(\epsilon)\subset N_{n}(\epsilon)$
\end{itemize}
We claim that for small $\epsilon_{0}>0$ whenever $\epsilon\le\epsilon_{0}$
and $n\ge n_{0}(\epsilon)$ then $N_{n}(\epsilon)$ is positive invariant
w.r.t. $f_{n}$. If this is not true then there is a sequence $\epsilon_{m}\to0$
and \[
y_{m}\in N_{n_{m}}(\epsilon_{m})\]
 with $f_{n_{m}}(y_{m})\notin N$. By definition of $N_{n_{m}}(\epsilon_{m})$
there is a sequence $\tilde{y}_{m}\in Y$, $x_{m}\in V(\epsilon_{m})$
and $k_{m}\ge0$ such that $d(y_{m},\tilde{y}_{m})<2^{-m}$, $f_{n_{m}}^{[0,k_{m}]}(x_{m})\subset\tilde{U}$
and $\tilde{y}_{m}=f_{n_{m}}^{k_{m}}(x_{m})$. Because $\phi(x_{m})\to0$
and $A_{f}^{-}(B)=A_{f}(B)$ we can assume w.l.o.g. that $x_{m}\to x_{0}\in A_{f}(B)$.
Admissibility and $f_{n_{m}}\to f$ imply the sequence $\{f_{n_{m}}^{k_{m}}(x_{m})\}_{m\in\mathbb{N}}$
has a convergent subsequence and w.l.o.g. $\tilde{y}_{m}=f_{n_{m}}^{k_{m}}(x_{m})\to y_{0}\in A_{f}^{-}(N')=A_{f}(N')\subset\operatorname{int}N$
and thus $y_{m}\to y_{0}$. 

Since $f_{n_{m}}(y_{m})\notin N\subset\operatorname{int}N'\cap f^{-1}(\operatorname{int}N')$
and $f_{n}\to f$ we have $f_{n_{m}}(y_{m})\to f(y_{0})\in N'\backslash\operatorname{int}N$.
But $y_{0}\in A_{f}(N')$ implies $f(y_{0})\in A_{f}(N')$ which is
a contradiction because $A_{f}(N')$ and $N'\backslash\operatorname{int}N$
are disjoint.\end{proof}
\begin{cor}
Under the assumption above for all $n\ge n_{0}$ we can find positive
$f_{n}$-invariant $N_{n},N_{n}^{'}$ such that \[
N_{n}\subset U_{\delta}(K)\subset N_{n}^{'}\subset B\]
for some $\delta$-neighborhood of $K$ denoted by \[
U_{\delta}(K)=\{x\in Y\,|\, d(x,y)<\delta\,\mbox{for some}\: y\in K\}.\]
Furthermore, we have $K\subset\operatorname{int}N_{n}$ and there
is an $\epsilon>0$ such that \[
U_{\epsilon}(K_{n})\subset N_{n}^{'}.\]
\end{cor}
\begin{proof}
Applying the previous theorem we get \[
K\cup K_{n}\subset V(\tilde{\epsilon})\subset N_{n}^{'}\subset B.\]

Recalling the definition of $\phi$ it is obvious that because $A_{f}^{-}(B)=A_{f}(B)=K$
for small $0<\epsilon'<1$ and $x\in P_{1}^{\epsilon'}$ \[
d(x,K)\le\epsilon'.\]
Because $K$ is compact and $V(\tilde{\epsilon})$ a neighborhood
of $K$ the $\delta$-neighborhood $U_{\delta}(K)$ of $K$ is contained
in $V(\tilde{\epsilon})$ for $\delta$ sufficiently small. Furthermore,
we can find an $\epsilon'>0$ with $\epsilon'<\delta$ such that $P_{1}^{\epsilon'}\subset U_{\delta}(K)\subset V(\tilde{\epsilon})$.
Applying the theorem again for $P_{1}^{\epsilon'}$ instead of $B$
we get for $n\ge n_{0}$ \[
N_{n}\subset P_{1}^{\epsilon'}\subset U_{\delta}(K)\subset N_{n}^{'}\subset B\]
and $K\subset V(\epsilon')\subset N_{n}$.

To show that $U_{\epsilon}(K_{n})\subset N_{n}$ we need another positive
$f_{n}$-invariant neighborhood $N_{n}^{''}$. First note that there
is a $\delta'>0$ such that \[
d(x,K)\ge\delta'\]
for all $x\in\partial P^{\epsilon'}$. So if we choose $0<2\epsilon<\delta$
then \[
U_{\epsilon}(P_{1}^{\epsilon})\subset P_{1}^{\epsilon'}.\]
Applying the previous theorem again we get a positive $f_{n}$-invariant
isolating $N_{n}^{''}$ of $K_{n}$ inside of $P_{1}^{\epsilon}$.
Hence \[
U_{\epsilon}(K_{n})\subset U_{\epsilon}(N_{n}^{''})\subset U_{\epsilon}(P_{1}^{\epsilon})\subset P_{1}^{\epsilon^{'}}\subset N_{n}^{'}.\]

\end{proof}
Now we are able to continue the attractor. Instead of an admissibility
assumption for the perturbed map we will use weak compactness of close
$\delta$-neighborhoods of compact sets.
\begin{defn}
[weak attractor]Suppose $Y$ has a weaker (Hausdorff) topology (i.e.
$x_{n}\to x$ strongly implies $x_{n}\rightharpoonup x$ weakly) and
$f$ is continuous and weakly continuous. An isolated invariant set
$K$ is called a weak attractor if it admits a positive $f$-invariant
isolating neighborhood $N$ such that $\omega^{\tiny\mbox{weak}}(N)\subset K$
where $\omega^{\tiny\mbox{weak}}(N)$ is defined as \[
\omega^{\tiny\mbox{weak}}(N)=\{y\in Y\,|\,\exists x_{n}\in N,m_{n}\to\infty\,\mbox{s.t. }\, f^{m_{n}}(x_{n})\rightharpoonup y\}.\]
\end{defn}
\begin{rem*}
(1) Our definition of weakness of an attractor is w.r.t. the weaker
topology and is different from one defined in \cite{Hurley2001}.
Even our definition of a (strong) attractor is weaker than the one
used there because we only require the existence of a positive invariant
isolating neighborhood of the invariant set. But there might be a
connection to Ochs' weak random attractor \cite{Ochs1999}.

(2) A Conley theory with weak-admissibility instead of admissibility
might not make sense since the continuation proof requires continuity
of the metric and usually the metric is only lower semicontinuous
w.r.t.  weak convergence.

(3) A weakly continuous function might not be continuous and vice
versa (see counterexample in the proof of theorem \ref{thm:wasserstein-cont})\end{rem*}
\begin{thm}
\label{thm:weak-continuation}Under the assumption of the previous
theorem, suppose there is a weaker (Hausdorff) topology on $Y$ and
that (strongly) closed $\delta$-neighborhoods of compact sets are
weakly (sequentially) compact, i.e. $\operatorname{cl}U_{\delta}(C)$
is weakly compact for compact $C$. If $f_{n}$ is weakly continuous
then $K_{n}$ is non-empty and a weakly compact weak attractor w.r.t.
$f_{n}$ for all $n\ge n_{0}$. Furthermore, $U_{\epsilon}(K_{n})\subset N_{n}^{'}$
for some $\epsilon>0$ and positive $f_{n}$-invariant $N_{n}^{'}$.\end{thm}
\begin{proof}
Applying the previous corollary we get \[
N_{n}\subset U_{\delta}(K)\subset N_{n}^{'}\subset B\]
and \[
U_{\epsilon}(K_{n})\subset N_{n}^{'}.\]
Because $f_{n}$ is weakly continuous, $N_{n}$ is positive $f_{n}$-invariant
and $\operatorname{cl}U_{\delta}(K)\subset N_{n}^{'}$ is closed and
thus weakly compact the set \[
\omega_{n}^{\tiny\mbox{weak}}(x)=\{y\in Y\,|\, f_{n}^{n_{k}}(x)\rightharpoonup y\,\mbox{for some}\, n_{k}\to\infty\}\]
 for $x\in N_{n}$ is non-empty and weakly compact. This implies $K_{n}\ne\varnothing$
and in particular \[
\omega_{n}^{\tiny\mbox{weak}}(x)\subset K_{n}\subset N_{n}.\]
Similarly weak compactness of $\operatorname{cl}U_{\delta}(K)$ implies
$\omega_{n}^{\tiny\mbox{weak}}(N_{n})\subset A_{n}(N_{n}^{'})=K_{n}\subset N_{n}$
and thus $K_{n}$ is a weak attractor. Obviously $K_{n}$ is weakly
closed and contained in the weakly compact set $\operatorname{cl}U_{\delta}(K)$
and is therefore weakly compact as well.
\end{proof}

\section{Wasserstein spaces}

Now we will introduce some notation and results for Wasserstein spaces
of a metric space, general references are \cite{AmbGigSav2008} and
\cite{Villani2009}.

Let $(X,d)$ be a complete separable metric space, also called Polish
space. We call it proper if every bounded closed set is compact. In
particular, this implies that $X$ is locally compact. The metric
of a non-compact proper metric space is necessarily unbounded.

The space of probability measures on the Borel $\sigma$-algebra of
$X$ is denoted by $\mathcal{P}(X)$. This space is given the weak
topology, i.e. $\mu_{n}\rightharpoonup\mu$ if $\int fd\mu_{n}\to\int fd\mu$
for all bounded continuous functions $f$ . Let $x_{0}$ be an arbitrary
point of $X$ and define $\mathcal{P}_{p}(X)$, the Wasserstein space
(of order $p$), by \[
\mathcal{P}_{p}(X)=\{\mu\in\mathcal{P}(X)\,|\,\int d(x_{0},x)^{p}d\mu(x)<\infty\}.\]
 Furthermore, define for $\mu,\nu\in\mathcal{P}_{p}(X)$\[
w_{p}(\mu,\nu)=\left(\inf_{\pi\in\Pi(\mu,\nu)}\int d(x,y)^{p}d\pi(x,y)\right)^{\frac{1}{p}}\]
where $\pi\in\Pi(\mu,\nu)\subset\mathcal{P}(X\times X)$ with $\pi(A\times X)=\mu$
and $\pi(X\times A)=\nu$ for all Borel sets $A$ and $B$. Then $(\mathcal{P}_{p}(X),w_{p})$
is a complete separable metric space. This topology is usually stronger
than the induced subspace topology of $\mathcal{P}_{p}(X)\subset\mathcal{P}(X)$. 

If $X$ is compact so is $\mathcal{P}_{p}(X)$. And $\mathcal{P}_{p}(X)$
is local compact only if $X$ is compact. A counterexample for non-proper
metric spaces is given in {\cite[7.1.9]{AmbGigSav2008}}. We will
adjust their example to non-compact proper metric spaces by showing
that the closed $\epsilon$-ball $B_{\epsilon}^{w_{p}}(\delta_{x_{0}})$
around $\delta_{x_{0}}$ in $\mathcal{P}_{p}(X)$ cannot be compact
for any $\epsilon>0$ and thus $\mathcal{P}_{p}(X)$ cannot be locally
compact. 
\begin{example*}
Assume $X$ is non-compact and proper and define \[
\mu_{n}=m_{n}\delta_{x_{n}}+(1-m_{n})\delta_{x_{0}}\]
 for some sequence $\{x_{n}\}_{x\in\mathbb{N}}\subset X$. Then \[
w_{p}(\mu_{n},\delta_{x_{0}})^{p}=m_{n}d(x_{n},x_{0})^{p}.\]
Suppose $d(x_{n},x_{0})\ge\epsilon>0$ and set $m_{n}=\epsilon^{p}\cdot d(x_{n},x_{0})^{-p}$
then $\mu_{n}\in\partial B_{\epsilon}^{w_{p}}(\delta_{x_{0}})$. If
$\{m_{n}\}_{n\in\mathbb{N}}$ stays bounded away from $0$ then $\{d(x_{n},x_{0})\}_{n\in\mathbb{N}}$
is bounded and thus $\{x_{n}\}_{n\in\mathbb{N}}$ and $\{m_{n}\}_{n\in\mathbb{N}}$
have convergent subsequences $x_{n'}\to x_{\infty}$ and $m_{n'}\to m\in(0,1]$
and $\mu_{n}\to m\delta_{x_{\infty}}+(1-m)\delta_{x_{0}}$ strongly
in $\mathcal{P}_{p}(X)$. But if we assume $d(x_{n},x_{0})\to\infty$
then $m_{n}\to0$ and thus $\mu_{n}\rightharpoonup\delta_{x_{0}}$
weakly. Because strong convergence requires that $w_{p}(\mu_{n},\delta_{x_{0}})=\epsilon$
converges to $0$ the sequence cannot converge strongly in $\mathcal{P}_{p}(X)$.
\end{example*}
Even though Wasserstein spaces are in general not locally compact
we can still show that the following holds for proper metric spaces.
The result is probably known or at least implicitly used in case $X=\mathbb{R}^{n}$.
Because it will be our main reason why the {}``weak'' Conley theory
is applicable and because we couldn't find any reference, we will
prove it completely.
\begin{thm}
\label{thm:weakly-proper}If $(X,d)$ is a proper metric space then
all closed $\delta$-neighborhoods of compact sets in $(\mathcal{P}_{1}(X),w_{1})$
are weakly compact, where the weak topology of $\mathcal{P}_{1}(X)$
is the induced subspace topology $\mathcal{P}_{1}(X)\subset\mathcal{P}(X)$.
A space, e.g. $\mathcal{P}_{1}(X)$, having this property may be called
weakly proper.\end{thm}
\begin{cor}
For $1\le p<q$ closed $\delta$-neighborhoods of compact sets in
$(\mathcal{P}_{q}(X),w_{q})$ are compact in $(\mathcal{P}_{p}(X),w_{p})$,
i.e. $\mathcal{P}_{q}(X)$ is weakly proper w.r.t. the subspace topology
induced by $\mathcal{P}_{q}(X)\subset\mathcal{P}_{p}(X)$. \end{cor}
\begin{rem*}
This is stronger then a compact embedding $i:(\mathcal{P}_{q}(X),w_{q})\to(\mathcal{P}_{p}(X),w_{p})$
because if $\mu_{n}\in\mathcal{P}_{q}(X)$ is bounded then w.l.o.g.
$i(\mu_{n})\to\mu_{*}$ in $\mathcal{P}_{p}(X)$ and necessarily $\mu_{*}\in\mathcal{P}_{q}(X)$,
i.e. bounded sequences never {}``leave'' the space.\end{rem*}
\begin{proof}
[Proof of theorem]We will show that \[
B_{r}^{w}:=\{\nu\in\mathcal{P}_{1}(X)\,|\, w_{1}(\nu,\delta_{x_{0}})\le r\}\]
 is weakly compact for all $r\ge0$. Since $B_{r}^{w}$ is closed
and $\mu_{n}\rightharpoonup\mu$ implies $w_{1}(\mu_{n},\delta_{x_{0}})\le\liminf_{n\to\infty}w_{1}(\mu_{n},\delta_{x_{0}})\le r$
we only need to show that $B_{r}^{w}$ is tight.

Tightness of a subset $\mathcal{K}\subset\mathcal{P}_{1}(X)$ means
for all $\epsilon>0$ there is a compact $K_{\epsilon}$ such that
for all $\mu\in\mathcal{K}$\[
\mu(X\backslash K_{\epsilon})\le\epsilon.\]

For $\mu\in B_{r}^{w}$ we have\[
\int_{X}d(x,x_{0})d\mu=w_{1}(\mu,\delta_{x_{0}})\le r.\]
Now choose $K_{\epsilon}=B_{\frac{r}{\epsilon}}(x_{0})$, the closed
ball around $x_{0}$ with radius $\frac{r}{\epsilon}$, which is compact
because $X$ is proper. Then we have \[
\mu(X\backslash K_{\epsilon})=\int_{X\backslash K_{\epsilon}}d\mu(x)\le\frac{\epsilon}{r}\int_{X\backslash K_{\epsilon}}d(x,x_{0})d\mu(x)\le\epsilon.\]
Thus the closed ball in $\mathcal{P}_{1}(X)$ around $\delta_{x_{0}}$
is weakly compact. 

Let $K\subset\mathcal{P}_{1}(X)$ be a compact set , e.g. $K=\{\mu\}$,
then the closed $R$-neighborhood around $K$ is defined as \[
N_{R}^{w}(K)=\{\nu\in\mathcal{P}_{1}(X)\,|\, w_{1}(\mu,\nu)\le R\,\mbox{ for some }\mu\in K\}.\]
 This set is closed and bounded and for some $\tilde{R}$ we have
\[
N_{R}^{w}(K)\subset B_{\tilde{R}}^{w}.\]

Let $\nu_{n}\in N_{R}^{w}(K)$ be an arbitrary sequence. Then there
are $\mu_{n}\in K$ with $w_{1}(\mu_{n},\nu_{n})\le R$. Because $B_{\tilde{R}}^{w}$
is weakly compact and $K$ is compact there are $\nu_{\infty}\in B_{\tilde{R}}^{w}$
and $\mu_{\infty}\in K$ such that for some subsequence (also denoted
by $\mu_{n}$, resp. $\nu_{n}$) \begin{eqnarray*}
\mu_{n} & \to & \mu_{\infty}\\
\nu_{n} & \rightharpoonup & \nu_{\infty}.\end{eqnarray*}
 Since $w_{1}(\cdot,\cdot)$ is weakly lower semicontinuity we have
\[
w_{1}(\mu,\nu)\le\liminf_{n\to\infty}w_{1}(\mu_{n},\nu_{n})\le R,\]
i.e. $\nu\in N_{R}^{w}(K)$ which implies weak compactness.
\end{proof}

\begin{proof}
[Proof of corollary] We only show that $B_{r}^{w_{q}}(\delta_{x_{0}})$
is weakly compact w.r.t. the induced subspace topology $\mathcal{P}_{q}(X)\subset\mathcal{P}_{p}(X)$
for $1\le p<q$. The rest will follow by the same arguments used above.

Assume $\{\mu_{n}\}_{n\in\mathbb{N}}\subset B_{r}^{w_{q}}(\delta_{x_{0}})$.
Since $w_{q}\le w_{1}$ the previous theorem implies w.l.o.g. $\mu_{n}\rightharpoonup\mu_{\infty}$
for some $\mu_{\infty}\in\mathcal{P}_{1}(X)$. Because \[
w_{q}(\mu_{\infty},\delta_{x_{0}})\le\liminf_{n\to\infty}w_{q}(\mu_{n},\delta_{x_{0}})\le r\]
we actually have $\mu_{\infty}\in B_{r}^{w_{q}}(\delta_{x_{0}})\subset\mathcal{P}_{q}(X)$. 

Because $1\le p<q$ \[
\int_{X\backslash B_{R}}d(x,x_{0})^{p}d\mu_{n}(x)\le\frac{1}{R^{q-p}}\int_{X\backslash B_{R}}d(x,x_{0})^{q}d\mu_{n}(x)\le\frac{r}{R^{q-p}}.\]
Hence \[
\lim_{R\to\infty}\limsup_{n\to\infty}\int_{X\backslash B_{R}}d(x,x_{0})^{p}d\mu_{n}(x)\le\lim_{R\to\infty}\frac{r}{R^{q-p}}=0.\]
This and the $\mu_{n}\rightharpoonup\mu_{\infty}$ weakly show that
$\mu_{n}\to\mu_{\infty}$ in $\mathcal{P}_{p}(X)$ (see {\cite[6.8]{Villani2009}}).
\end{proof}
In the following assume that $X$ is proper. Suppose now $f:X\to X$
is a continuous map having a global (set) attractor, i.e. there is
a compact $f$-invariant $A\subset X$ such that for all bounded sets
$B$\[
\lim_{n\to\infty}\tilde{d}(f^{n}(B),A)=0,\]
where $\tilde{d}$ is the semi-Hausdorff metric induced by $d$ such
that $\tilde{d}(A,B)=0$ iff $A\subset\operatorname{cl}B$. The map
$f$ induces a continuous map $f_{*}:\mathcal{P}(X)\to\mathcal{P}(X)$
with $f_{*}(\mu)(B)=\mu(f^{-1}(B))$ for all Borel set $B$. Furthermore,
under slightly stronger assumptions $f_{*}$ has a global attractor
\[
\mathcal{P}(A)=\{\mu\in\mathcal{P}(X)\,|\,\mu(A)=1\}.\]

\begin{rem*}
For compact $X$ the global attractor is always $X$ itself. In particular,
since $\mathcal{P}_{p}(X)=\mathcal{P}(X)$ for $1\le p<\infty$ is
compact the global attractor of $\mathcal{P}(X)$ is the space itself
and we don't get new information. The whole theory is only interesting
for non-compact proper metric spaces $X$.
\end{rem*}
Since the Wasserstein space includes distance the behavior of $f$
at infinity becomes important. 
\begin{thm}
\label{thm:wasserstein-cont}The map $f_{*}$ induces a continuous
map $(\mathcal{P}_{p}(X),w_{p})\to(\mathcal{P}_{p}(X),w_{p})$ (also
denoted by $f_{*}$) if and only if for some $x_{0}\in X$ \[
\sup_{x\in X}\frac{d(f(x),x_{0})}{1+d(x,x_{0})}<\infty.\]
\end{thm}
\begin{rem*}
For a semiflow $\pi$ we need that $\sup_{x\in X}\nicefrac{d(x\pi t,x_{0})}{1+d(x,x_{0})}=M_{t}<\infty$
for $t\in[0,T_{0}]$. Which means, in particular, that there has to
be a global lower bound on the blow-up time and thus there cannot
be blow-ups at all, i.e. $\pi$ has to be a global semiflow. Which
implies that the induced semiflow $\pi_{*}$ on $\mathcal{P}_{p}(X)$
is also a global semiflow. A necessary requirement for the existence
of a global attractor is an upper bound on $M_{t}$ for all $t\ge0$.
The requirements in {\cite[Chapter 8]{AmbGigSav2008}} are sometimes
too strong. A sufficient condition is that a one-sided Lipschitz condition
holds globally (e.g. $v(x)=x-x^{3}$ is an unbounded vector field
and only locally Lipschitz, but satisfies a one-sided Lipschitz condition).\end{rem*}
\begin{proof}
Suppose first that \[
M=\sup_{x\in X}\frac{d(f(x),x_{0})}{1+d(x,x_{0})}<\infty.\]
We can assume w.l.o.g. $M>0$, otherwise $f|_{X}\equiv x_{0}$ and
the result is obvious. For $\mu\in\mathcal{P}_{p}(X)$ we have \begin{eqnarray*}
\int d(x,x_{0})^{p}df_{*}\mu(x) & = & \int d(f(x),x_{0})^{p}d\mu(x)\\
 & \le & M^{p}\int(1+d(x,x_{0}))^{p}d\mu(x)<\infty,\end{eqnarray*}
 i.e. $f_{*}(\mu)\in\mathcal{P}_{p}(X)$. So we only need to show
continuity.

Suppose $\mu_{n}\to\mu$ in $\mathcal{P}_{p}(X)$ then $f_{*}\mu_{n}\rightharpoonup f_{*}\mu$.
Since $d(f(\cdot),x_{0})^{p}$ is continuous and grows at most like
$d(\cdot,x_{0})^{p}$ it follows that\begin{eqnarray*}
\int d(x,x_{0})^{p}df_{*}\mu_{n}(x) & = & \int d(f(x),x_{0})^{p}d\mu_{n}(x)\\
 &  & \longrightarrow\int d(f(x),x_{0})^{p}d\mu(x)=\int d(x,x_{0})^{p}df_{*}\mu(x).\end{eqnarray*}
 Thus $f_{*}\mu_{n}\to f_{*}\mu$ (see {\cite[6.8]{Villani2009}})
which shows that $f_{*}:(\mathcal{P}_{p}(X),w_{p})\to(\mathcal{P}_{p}(X),w_{p})$
is (strongly) continuous. 

It remains to show that $f_{*}$ is not continuous if there is a sequence
$\{x_{n}\}_{n\in\mathbb{N}}\subset X$ such that \[
d_{n}=\frac{d(f(x_{n}),x_{0})}{1+d(x_{n},x_{0})}\to\infty.\]
 Because $X$ is proper and $f$ continuous we must have $d(x_{n},x_{0})\to\infty$.
For large $n$ we can assume $0<\frac{1}{d_{n}}<d(x_{n},x_{0})$.
Set $c_{n}=\frac{1}{d_{n}}$ then \[
\mu_{n}=c_{n}^{p}\frac{1}{d(x_{n},x_{0})^{p}}\delta_{x_{n}}+(1-c_{n}^{p}\frac{1}{d(x_{n},x_{0})^{p}})\delta_{x_{0}}\in\mathcal{P}_{1}(X)\]
and $w_{p}(\mu_{n},\delta_{x_{0}})^{p}=c_{n}^{p}\to0$, i.e. $\mu_{n}\to\delta_{x_{0}}$
strongly in $\mathcal{P}_{p}(X)$. We have \[
f_{*}\mu_{n}=c_{n}^{p}\frac{1}{d(x_{n},x_{0})^{p}}\delta_{f(x_{n})}+(1-c_{n}^{p}\frac{1}{d(x_{n},x_{0})^{p}})\delta_{f(x_{0})}\]
and therefore \[
w_{p}(f_{*}\mu_{n},f_{*}\delta_{x_{0}})^{p}=c_{n}^{p}\frac{d(f(x_{n}),x_{0})^{p}}{d(x_{n},x_{0})^{p}}=\frac{1+d(x_{n},x_{0})^{p}}{d(x_{n},x_{0})^{p}}\to1\]
which implies that $f_{*}$ cannot be (strongly) continuous in $\delta_{x_{0}}$.
\end{proof}
The following results will hold for any $\mathcal{P}_{p}(X)$. To
simplify the notation and some of the proofs we will just state them
for $\mathcal{P}_{1}(X)$. Furthermore, we assume from now on that
$f_{*}:\mathcal{P}_{1}(X)\to\mathcal{P}_{1}(X)$ is strongly continuous
(for short just continuous) and whenever we speak about $f_{*}$ we
mean the map $f_{*}:\mathcal{P}_{1}(X)\to\mathcal{P}_{1}(X)$. Since
$f_{*}(\mathcal{P}_{1}(X))\subset\mathcal{P}_{1}(X)$ (for $f_{*}:\mathcal{P}(X)\to\mathcal{P}(X)$)
this also implies that $f_{*}$ is weakly continuous in $\mathcal{P}_{1}(X)$.
Similarly we could say that $f_{*}:\mathcal{P}_{q}(X)\to\mathcal{P}_{q}(X)$
is continuous and {}``weakly'' continuous in $\mathcal{P}_{q}(X)$
w.r.t. the induced subspace topology of $\mathcal{P}_{q}(X)\subset\mathcal{P}_{p}(X)$
for any $1\le p<q$.
\begin{example*}
Having a global attractor does not imply that $f_{*}$ is strongly
continuous, even finite time compactness is not sufficient: Let $X$
be $\mathbb{R}$ with the Euclidean metric $|\cdot|$. Define $f:\mathbb{R}\to\mathbb{R}$
by \[
f(x)=\begin{cases}
0 & x\ge0\\
x^{2} & x<0.\end{cases}\]
Then $f$ is continuous and $f^{2}\equiv0$ but for $x_{n}=-n$ \[
\frac{d(f(x_{n}),0)}{1+d(x_{n},0)}=\frac{n^{2}}{1+n}\to\infty,\]
i.e. $f_{*}$ is not continuous on $(\mathcal{P}_{p}(X),w_{p})$.
\end{example*}
If $K\subset\mathcal{P}_{1}(X)$ is invariant w.r.t. $f_{*}$ then
it is invariant w.r.t. $f_{*}:\mathcal{P}(X)\to\mathcal{P}(X)$. Which
implies that all measures in $K$ are supported on the global attractor
$A$ of $f$. Since $\mathcal{P}(A)\subset\mathcal{P}_{1}(X)$ the
maximal invariant set of $\mathcal{P}_{1}(X)$ is $\mathcal{P}(A)=\mathcal{P}_{1}(A)$.

Suppose $f$ is finite time compact, i.e. there is an $m$ such that
$f^{m}(X)\subset B_{R}$ for some compact set $B_{R}$. It should
be obvious that this implies $K=\mathcal{P}_{1}(A)$ is the global
attractor of $f_{*}$. Furthermore, we have the following:
\begin{prop}
Suppose for some $m>0$, \[
f^{m}(X)\subset B_{R}(x_{0}).\]
Then $f_{*}$ is finite time compact and thus any closed set $B^{w}\subset\mathcal{P}_{1}(X)$
is $f_{*}$-admissible. \end{prop}
\begin{proof}
Let $\{\nu_{n}\}_{n\in\mathbb{N}}$ be any sequence in $f_{*}^{m}(\mathcal{P}_{1}(X))$.
Then there is a sequence $\{\mu_{n}\}_{n\in\mathbb{N}}$ such that
$\nu_{n}=f_{*}^{m}(\mu_{n})$. $f^{m}(X)\subset B_{R}(x_{0})$ implies
$\operatorname{supp}\nu_{n}\subset B_{R}(x_{0})$. Thus $\nu_{n}$
is tight and \[
\int_{d(x,x_{0})\ge R+\epsilon}d(x,x_{0})d\nu_{n}(x)=0,\]
i.e. $\{\nu_{n}\}_{n\in\mathbb{N}}$ has uniformly integrable first
moments. Which means that $\{\nu_{n}\}_{n\in\mathbb{N}}$ has a convergent
subsequence. Therefore, $f_{*}^{m}(\mathcal{P}_{1}(X))$ is compact
which easily implies admissibility for any closed $B^{w}$.\end{proof}
\begin{prop}
Suppose there is an $R_{0}$, $0\le c<1$ and $m>0$ such that for
all $R\ge R_{0}$\[
f^{m}(B_{R})\subset B_{cR}.\]
Then any bounded closed set $B^{w}\subset\mathcal{P}_{1}(X)$ is $f_{*}$-admissible.\end{prop}
\begin{proof}
Let $\{\mu_{n}\}_{n\in\mathbb{N}}$ be a sequence in $B^{w}$ such
that $f_{*}^{[0,m_{n}]}(\mu_{n})\subset B^{w}$ for some $m_{n}\to\infty$.
Because $B^{w}$ is bounded we have\[
\int_{X}d(x,x_{0})d\mu_{n}(x)\le M.\]

First assume $m_{n}=k_{n}\cdot m$ for an unbounded sequence $k_{n}\in\mathbb{N}$.
Then for $R\ge R_{0}$ \begin{eqnarray*}
\int\chi_{X\backslash B_{R}}(x)\cdot d(x,x_{0})df_{*}^{m_{n}}\mu_{n}(x) & = & \int\chi_{X\backslash B_{R}}(f^{m_{n}}(x))\cdot d(f^{m_{n}}(x),x_{0})d\mu_{n}(x)\\
 & \le & \int\chi_{X\backslash B_{c^{-k_{n}}R}}(x)\cdot c^{k_{n}}d(x,x_{0})d\mu_{n}(x)\\
 & \le & c^{k_{n}}M\to0,\end{eqnarray*}
which shows that $\{f_{*}^{m_{n}}(\mu_{n})\}_{n\in\mathbb{N}}$ has
uniformly integrable first moments which implies that the sequence
has a convergent subsequence. 

If $m_{n}\not\equiv0(\operatorname{mod}m)$ then for some $0\le l_{n}<m$
we have $m_{n}-l_{n}\equiv0(\operatorname{mod}m)$. Therefore, if
we set $\nu_{n}=f_{*}^{l_{n}}(\mu)$ then the argument above applies
to $\nu_{n}$ and the sequence of endpoints (which is equal to $\{f_{*}^{m_{n}}(\mu_{n})\}_{n\in\mathbb{N}}$)
has a convergent subsequence.\end{proof}
\begin{rem*}
Kifer used in {\cite[Theorem 1.7]{Kifer1988}} linear attraction instead
of exponential. This might not be sufficient for admissibility. Nevertheless,
later we will assume that a Markov-type perturbation of $f_{*}$ is
small in the Wasserstein distance which is stronger than Kifer's assumption
and thus an invariant (probability) measure exists for the perturbation
by the same theorem. But that theorem does not imply that the perturbed
and unperturbed invariant measures are close w.r.t. the Wasserstein
distance, the perturbed invariant measures might not even be in the
Wasserstein space. So our result improves this sufficiently.
\end{rem*}
Before we show how to use Conley theory for small Markov-type noise
applied to $f$ we give a sufficient condition such that bounded sets
in $\mathcal{P}_{1}(X)$ are $\{F_{n}\}$-admissible for $F_{n}\to f_{*}$.
\begin{prop}
Let $B^{w}$ be closed and bounded and $U$ be a $\delta$-neighborhood
of $B^{w}$ with $f_{*}$-admissible closure. Suppose $F_{n}\to f_{*}$
uniformly on some $U$, i.e. \[
\sup_{\mu\in U}w_{1}(f_{*}(\mu),F_{n}(\mu))=\epsilon_{n}\to0.\]
If $f_{*}$ is uniformly continuous in $U$ then $B^{w}$ is $\{F_{n}\}$-admissible.\end{prop}
\begin{rem*}
(1) The idea is to use the uniform convergence and uniform continuity
to construct longer and longer orbits of $f_{*}$ close the the last
part of the orbits of $F_{n}$, i.e. the orbit $f^{[0,m_{n}-k_{n}]}(y_{n})$
and $F_{n}^{[k_{n},m_{n}]}(x_{n})$ should be closer and closer and
$m_{n}-k_{n}\to\infty$ for $y_{n}=F_{n}^{k_{n}}(x_{n})$.

(2) Uniform continuity of $f_{*}$ and uniform convergence of $F_{n}\to f_{*}$
are the assumptions Benci \cite{Benci1991} used to prove his continuation
theorem for the Conley index. Besides having continuous time dynamical
systems he also needs invertibility. \end{rem*}
\begin{proof}
Let $\mu_{n}\in N$ and $m_{n}\to\infty$ be sequences with $F_{n}^{[0,m_{n}]}(\mu_{n})$.
Uniform convergence of $F_{n}\to f_{*}$ and uniform continuity of
$f$ imply that for some $\epsilon(\epsilon_{n})\to0$ as $\epsilon_{n}\to0$
\begin{eqnarray*}
w_{1}(F_{n}^{2}(\mu),f_{*}^{2}(\mu)) & \le & w_{1}(F_{n}^{2}(\mu),f_{*}(F_{n}(\mu))+w_{1}(f_{*}(F_{n}(\mu)),f_{*}^{2}(\mu))\\
 & \le & \epsilon_{n}+\epsilon(\epsilon_{n})=:\epsilon_{n,2}\to0\,\mbox{as \,}n\to\infty.\end{eqnarray*}
Similarly we can show that there are $\epsilon_{n,k}\to0$ as $n\to\infty$
such that \[
w_{1}(F_{n}^{k}(\mu),f_{*}^{k}(\mu))\le\epsilon_{n,k}\to0.\]

Therefore, there is a sequence $k_{n}\ge0$ with $m_{n}-k_{n}\to\infty$
such that $F_{n}^{[0,m_{n}]}(\mu_{n})\subset B^{w}$ implies that
\[
f_{*}^{[0,m_{n}-k_{n}]}(\nu_{n})\subset U=U_{\delta}(B^{w})\]
for $\nu_{n}=F_{n}^{k_{n}}(\mu_{n})$. Furthermore, we can choose
$k_{n}$ such that \[
\delta_{n}=\max_{k\in[0,m_{n}-k_{n}]}\epsilon_{n,k}\to0\]
 and therefore \[
w_{1}(F_{n}^{m_{n}}(\mu_{n}),f_{*}^{m_{n}-k_{n}}(\nu_{n}))\le\delta_{n}.\]

Because the closure of $U$ is $f_{*}$-admissible, the sequence of
endpoints $\{f_{*}^{m_{n}-k_{n}}(\nu_{n})\}_{n\in\mathbb{N}}\subset U$
has a convergent subsequence which implies that $\{F_{n}^{m_{n}}(\mu_{n})\}_{n\in\mathbb{N}}$
has a convergent subsequence, in particular the limit point is in
$B^{w}$.
\end{proof}
This proposition applies in particular to Lipschitz continuous functions
$f:X\to X$ because the induce map $f_{*}:\mathcal{P}_{1}(X)\to\mathcal{P}_{1}(X)$
is Lipschitz continuous as well.

E.g. suppose $X=\mathbb{R}^{n}$ and $f$ is the time $h$ map of
a flow generated by an ODE $\dot{x}=g(x)$ such that $g$ satisfies
the one-sided Lipschitz condition for some $M\in\mathbb{R}$\[
\langle x-y,g(x)-g(y)\rangle\le M\|x-y\|^{2}\]
then $f$ is Lipschitz continuous with constant $e^{Mh}$.

\section{\label{sec:MW-maps}Markov-Wasserstein maps}
\begin{defn}
A Markov-Wasserstein map (MW-map) of order $p$ is a continuous map
$P:\mathcal{P}_{p}(X)\to\mathcal{P}_{p}(X)$ which is convex linear,
i.e. for $\mu,\nu\in\mathcal{P}_{p}(X)$ and $a\in[0,1]$\[
P(a\mu+(1-a)\nu)=aP(\mu)+(1-a)P(\nu).\]
Suppose $P$ is induced by a kernel $p(dy|x):X\to\mathcal{P}_{p}(X)$
(necessarily continuous), i.e. \[
P:d\mu(y)\mapsto\int p(dy|x)d\mu(x).\]
 The map $P=P_{f}^{M}$ is called an MW-map of order $p$ relative
to $f$ with noise level (at most) $M$ if\[
\sup_{x\in X}w_{p}(p(dy|x),\delta_{f(x)})\le M.\]
This implies by {\cite[4.8]{Villani2009}} \[
w_{p}(P(\mu),f_{*}(\mu))^{p}\le\int w_{p}(p(dy|x),\delta_{f(x)})^{p}d\mu(x)\le M^{p}.\]
\end{defn}
\begin{rem*}
As in the sections before, the results also hold for semiflows and
a suitable definition for MW-semiflows, i.e. a continuous semigroups
$(P_{t})_{t\ge0}$ on $\mathcal{P}_{p}(X)$. The noise level model
can be stated similarly, but we only require that it is uniformly
small for all {}``small'' $t$. We will focus here only on maps,
resp. MW-maps, because the intuition behind these is easier, but all
results also hold for MW-semiflows if the noise level is sufficiently
small.
\end{rem*}
MW-maps (resp. MW-semiflows) appear naturally in the theory of Markov
chains (resp. processes). The map $p(dy|\cdot):X\to\mathcal{P}_{p}(X)$
is the Markov transition probability function, whereas the map $P:\mathcal{P}_{p}(X)\to\mathcal{P}_{p}(X)$
almost never has a name. If $P=f_{*}$ for some dynamical system $f:X\to X$
then $P$ is sometimes called transfer map. The following will show
that we only need continuity of $p(dy|\cdot)$ to ensure that $P$
is continuous in $\mathcal{P}_{p}(X)$ (and thus for any $\mathcal{P}_{q}(X)$,
$1\le q<p$ and for $\mathcal{P}(X)$).
\begin{thm}
Let $p(dy|\cdot):X\to\mathcal{P}_{p}(X)$ be a Markov kernel, i.e.
a measure-valued map, continuously depending on $x$. If $M=\sup_{x\in X}w_{p}(p(dy|x),\delta_{x})<\infty$
then $P$ defined by \[
P:d\mu(y)\mapsto\int p(dy|x)d\mu(x)\]
 is an MW-map of order $p$ relative to $\operatorname{id}:X\to X$
with noise level $M$.\end{thm}
\begin{rem*}
For MW-semiflows weak continuity of $p_{t}(dy|\cdot)$ corresponds
to Feller continuity of the corresponding stochastic process. In fact,
if the initial distribution of $(X_{t}^{n})_{t\ge0}$ is $\delta_{x_{n}}$,
i.e. $X_{0}^{n}=x_{n}$, then by our continuity requirement if $x_{n}\to x_{0}$
then $\mu_{t}^{n}=P_{t}\delta_{x_{t}}\to P_{t}\delta_{x_{0}}=\mu_{t}^{0}$
which implies that \[
u(x_{n})=\mathbb{E}g(X_{t}^{n})=\int g(x)d\mu_{t}^{n}(x)\to\int g(x)d\mu_{t}^{0}(x)=\mathbb{E}g(X_{t}^{0})=u(x_{0})\]
 for all bounded continuous function $g(x)$, i.e. the stochastic
process generated by $(P_{t})_{t\ge0}$ is Feller continuous. The
continuity of the moments implies that, in addition, the moments are
also continuous. This condition could be called $p$-Feller continuous.
This type of continuity does not need $t\in\mathbb{R}$ and thus applies
equally to Markov chains, i.e. discrete time stochastic processes.\end{rem*}
\begin{proof}
Continuous dependency implies that $P:\mathcal{P}_{p}(X)\to\mathcal{P}(X)$
is continuous. So we only need to show that $\mu_{n}\to\mu$ in $\mathcal{P}_{p}(X)$
implies that $P\mu_{n}\to P\mu$ in $\mathcal{P}_{p}(X)$

Because $w_{p}(p(dy|x),\delta_{x})\le M$ for some $M<\infty$ we
have\begin{eqnarray*}
\int d(y,x_{0})^{p}p(dy|x) & = & w_{p}(p(dy|x),\delta_{x_{0}})^{p}\\
 & \le & (w_{p}(p(dy|x),\delta_{x})+w_{p}(\delta_{x},\delta_{x_{0}}))^{p}\\
 & \le & 2^{p-1}(M^{p}+d(x,x_{0})^{p})\end{eqnarray*}

This implies that $g(x)=\int d(y,x_{0})^{p}p(dy|x)$ grows at most
like $d(x,x_{0})^{p}$. Furthermore, $g$ is continuous because $x\mapsto p(dy|x)$
and $\mu\mapsto w_{p}(\mu,\delta_{x_{0}})^{p}$ are. Thus by {\cite[6.8]{Villani2009}}
\[
\int d(y,x_{0})^{p}dP\mu_{n}(x)=\int g(x)d\mu_{n}(x)\to\int g(x)d\mu(x)=\int d(y,x_{0})^{p}dP\mu(y),\]
i.e. $P\mu_{n}\to P\mu$ in $\mathcal{P}_{p}(X)$.\end{proof}
\begin{example*}
(1) Bounded noise can be modeled by Markov maps with \[
M=\sup_{x\in X}d(x,\operatorname{supp}p(dy|x))<\infty.\]
 Then $w_{p}(p(dy|x),\delta_{x})\le M$ and thus continuity of $p(dy|\cdot):X\to\mathcal{P}_{p}(X)$
for some $p$ implies that of $P:\mathcal{P}_{p}(X)\to\mathcal{P}_{p}(X)$.
This could also be used to model multi-valued perturbations, i.e.
maps $f_{n}:X\to2^{X}$ with $\sup_{x\in X}d(f(x),f_{n}(x))\le M$. 

(2) Let $X=\mathbb{R}^{n}$ with its Euclidean distance. If $\nu$
is the standard normal distribution then $\nu=\rho(x)dx$, where $dx$
is the Lebesgue measure on $\mathbb{R}^{n}$. Any normal distribution
with mean $x$ and variance $\sigma^{2}$ can be modeled as follows
\[
\mu_{x,\sigma^{2}}=\delta_{x}*\rho_{\sigma}\]
where $\rho_{\epsilon}(x)=\epsilon^{-n}\rho(x/\epsilon)$. For $\sigma=0$
we set $\mu_{x,0}=\delta_{x}$. 

Since $m_{p}=\int|x|^{p}\rho(x)dx<\infty$ for all $p$ this implies
(see {\cite[7.1.10]{AmbGigSav2008}}) that \[
w_{p}(\delta_{x},\mu_{x,\sigma^{2}})\le\sigma m_{p}.\]
Thus Gaussian noise with uniformly small variance is uniformly small
in all Wasserstein spaces (although the noise level diverges to $\infty$
as $p\to\infty$).\end{example*}
\begin{cor}
If $f:X\to X$ induces a continuous self map on $\mathcal{P}_{p}(X)$
and $p(dy|\cdot):X\to\mathcal{P}_{p}(X)$ is continuous with \[
M=\sup_{x\in X}w_{p}(p(dy|x),\delta_{f(x)})<\infty\]
 then $P:\mathcal{P}_{p}(X)\to\mathcal{P}_{p}(X)$ defined as above
is an MW map of order $p$ relative to $f$ with noise level $M$.
\end{cor}
A random perturbation can now be modeled as a composition of $f$
followed by a smearing via $p(dy|\cdot)$, i.e. $f_{*}$ followed
by $P=P_{\operatorname{id}}^{M}$. This corresponds to additive noise
depending only on the image, whereas a general MW-map relative to
$f$ might smear the image $f(x)$ and $f(y)$ for $x\ne y$ differently
even if $f(x)=f(y)$.

If for some sequence $\tilde{p}_{n}(dy|\cdot)$ the noise level\[
\sup_{x\in X}w_{p}(\tilde{p}_{n}(dy|x),\delta_{f(x)})=\epsilon_{n}\]
converges to zero then $\tilde{P}_{n}$ converges to $f_{*}$ uniformly
on $\mathcal{P}_{p}(X)$, i.e.\[
\sup_{\mu\in\mathcal{P}_{p}(X)}w_{p}(\tilde{P}_{n}(\mu),f_{*}(\mu))\to0.\]

An MW-chain relative to $f$ satisfies the Markov property, i.e. future
behavior only depends on the current state. Furthermore, this models
only time-independent random perturbations. Time-dependent perturbations
can be modeled with the result of \cite{Kell2011b}. There it is shown
that a local attractor can be continued if the non-autonomous perturbations
is uniformly small. Translated into this framework this means \[
\sup_{\mu\in\mathcal{P}_{p}(X),k\in\mathbb{Z}}w_{p}(P(\mu,k),f_{*}(\mu))<\epsilon\]
for the non-autonomous dynamical system ($\approx$ inhomogeneous
Markov map)\[
(\mu,k)\mapsto(P(\mu,k),k+1).\]
Instead of using the semi-admissibility argument to show that the
invariant set $K_{n}$ is non-empty we can use a weak compactness
argument to get the same result.
\begin{example*}
Consider the ODE with $\dot{x}=x-x^{3}$. This satisfies the one-side
Lipschitz condition with $M=1$ and generates a global semiflow that
attracts in finite time. Thus the time-one map induces a Lipschitz
continuous map $f_{*}$ on $\mathcal{P}_{p}(X)$ which attracts in
finite time, too. Therefore, any bounded closed set in $\mathcal{P}_{p}(X)$
is $\{F_{n}\}$-admissible for $F_{n}\to f_{*}$ uniformly. In particular,
the MW-map of order $p$ for small noise level has an attractor close
to the original w.r.t. the Wasserstein metric $w_{p}$. \end{example*}
\begin{thm}
\label{thm:weak-cont-wasserstein}Suppose $f:X\to X$ induces a dynamical
system $f_{*}$ on $\mathcal{P}_{p}(X)$ having a global attractor
and that $f_{*}$ is uniformly continuous in a neighborhood of the
global attractor. If $P_{n}\to f_{*}$ is a sequence of MW-maps of
order $p$ relative to $f$ with noise level $\epsilon_{n}\to0$.
Then for $n\ge n_{0}$ there is a positive $P_{n}$-invariant isolating
neighborhood $N_{n}$ such that $K_{n}=A_{P_{n}}(N_{n})$ is non-empty
and a weakly compact weak attractor which contains all bounded $P_{n}$-invariant
measures, i.e. $A_{P_{n}}(P_{p}(X))=K_{n}$. Furthermore, there is
at least one stationary measure in $K_{n}$.\end{thm}
\begin{rem*}
Suppose $p>1$. Whenever $K$ is (strongly) compact in $\mathcal{P}_{p}(X)$
then $\operatorname{cl}U_{\delta}(K)$ is weakly compact w.r.t. the
weaker subspace topology of $\mathcal{P}_{p}(X)$ induced by $\mathcal{P}_{p}(X)\subset\mathcal{P}_{q}(X)$
for any $1\le q<p$. Thus for all $\mu\in N_{n}$ there is a $\mu_{K}\in K_{n}$
\[
P_{n}(\mu)\overset{w_{q}}{\longrightarrow}\mu_{K}.\]
This means that, although the $p$-moment may not converge, any $q$-moment
converges for $1\le q<p$, but the convergence may get worse the closer
$q$ comes to $p$.\end{rem*}
\begin{proof}
Everything but the existence of a stationary measure and $A_{P_{n}}(\mathcal{P}_{p}(X))=K_{n}$
follows from theorem \ref{thm:weak-continuation}. Since $P_{n}$
is convex linear $K_{n}$ must be convex. This implies that for any
$\mu\in K_{n}$ the sequence \[
\left\{ \frac{1}{m}\sum_{k=0}^{m-1}P_{n}^{k}(\mu)\right\} _{m\in\mathbb{N}}\]
is in $K_{n}$ and thus weakly converging to some $\nu\in K_{n}$
and by the Krylov-Bogolyubov theorem it must be a fixed point of $P_{n}$,
i.e. $\nu$ is a stationary measure of $P_{n}$. Furthermore, if $R$
is the distance from the global attractor of $f$ then its mass must
decay as $R^{-p}$. 

The invariant measures must all be contained in the interior of $N_{n}$.
Otherwise take $\mu\in A_{P_{n}}(\mathcal{P}_{p}(X))\backslash K_{n}$.
If $\mu$ is $P_{n}$-stationary then the argument is as follows:
For $t\in[0,1]$ and some $\mu_{0}\in K_{n}$ the graph of $t\mapsto t\mu+(1-t)\mu_{0}$
is stationary and intersects $\partial N_{n}$, which implies $K_{n}=A_{P_{n}}(N_{n})$
intersects $\partial N_{n}$. This contradicts the isolatedness of
$K_{n}$.

For the general case assume $\sigma:\mathbb{Z}\to\mathcal{P}_{p}(X)$
is a bounded full solution through $\mu$, i.e. $P_{n}(\sigma(k))=\sigma(k+1)$
and $\sigma(0)=\mu$. Now define the function $g:\mathbb{Z}\to[0,1]$
\[
g:k\mapsto\sup\{t\in[0,1]\,|\, s\sigma(k)+(1-s)\mu_{0}\in N_{n}^{'}\,\mbox{for all}\, s\in[0,t]\}.\]
Because $N_{n}^{'}$ is a positive invariant neighborhood of $K_{n}$,
$\mu_{0}$ stationary, $P_{n}$ convex linear and $\{\sigma(k)\}_{k\in\mathbb{Z}}$
bounded and not entirely in $N_{n}^{'}$ we have \[
0<\delta\le g(k)\le g(k+1).\]
This implies that \[
T=\inf_{k\le0}g(k)\ge\delta.\]
Because $\sigma(0)\notin N_{n}^{'}$ we have $T<1$ and thus by definition
of $g$ \[
w_{p}(T\sigma(k)+(1-T)\mu_{0},\partial N_{n}^{'})\to0\quad\mbox{as\,}k\to-\infty\]
and thus there is a $T_{1}\le T$ such that $\tilde{\sigma}(k)=T_{1}\sigma(k)+(1-T_{1})\mu_{0}$
is a full solution in $N_{n}^{'}$ with \[
w_{p}(\tilde{\sigma}(k_{1}),\partial N_{n}^{'})\le\frac{\epsilon}{2}\]
for some $k_{1}\le0$. But $\tilde{\sigma}(k_{1})\in K_{n}$ and $U_{\epsilon}(K_{n})\subset N_{n}^{'}$
which implies that \[
\epsilon\le w_{p}(\sigma(k_{1}),\partial N_{n}^{'})\le\frac{\epsilon}{2}.\]
This is a contradiction and thus $A_{P_{n}}(\mathcal{P}_{p}(X))\backslash K_{n}=\varnothing$,
i.e. $K_{n}$ contains all bounded invariant measures in $\mathcal{P}_{p}(X)$.\end{proof}
\begin{cor}
The positive $f_{*}$-invariant isolating neighborhood $B$ and the
positive $P_{n}$-invariant isolating neighborhood $N_{n}$ can be
chosen convex, i.e. if $\mu_{i}\in B$ (resp. $\mu_{i}\in N_{n}$)
for $i=0,1$ then $\mu_{t}\in B$ (resp. $\mu_{t}\in N_{n}$) for
$\mu_{t}=t\mu_{0}+(1-t)\mu_{1}$ and $t\in[0,1]$.\end{cor}
\begin{proof}
Let $\nu_{i},\mu_{i}\in\mathcal{P}_{p}(X)$ for $i=0,1$ and define
$\mu_{t}=t\mu_{1}+(1-t)\mu_{0}$ and $\nu_{t}=t\nu_{1}+(1-t)\nu_{0}$.
Assume \[
w_{p}(\mu_{i},\nu_{i})<\epsilon.\]
Then there are optimal transference plans $\pi_{i}\in\Pi(\mu_{i},\nu)$
such that\[
\int d(x,y)^{p}d\pi_{i}(x,y)<\epsilon^{p}.\]
The plan $\pi_{t}=t\pi_{1}+(1-t)\pi_{0}$ is a transference plan for
the pair $(\mu_{t},\nu_{t})$ and thus\begin{eqnarray*}
w_{p}(\mu_{t},\nu_{t})^{p} & \le & \int d(x,y)^{p}d\pi_{t}(x,y)\\
 & \le & t\int d(x,y)^{p}d\pi_{1}(x,y)+(1-t)\int d(x,y)^{p}d\pi_{0}(x,y)\\
 & < & t\epsilon^{p}+(1-t)\epsilon^{p}=\epsilon^{p}.\end{eqnarray*}
The construction of $B$ is done via a Lyapunov pair $(\phi,\gamma)$
essentially measuring a weighted distance of the forward orbit of
a point, i.e. \[
F_{N'}:\mu\mapsto\min\{1,w_{p}(\mu,A^{-}(N')\cup\partial N')\}\]
and\[
\phi:\mu\mapsto\sup\{(2n+1)F_{N'}(f_{*}^{n}(x))/(n+1)\,|\, n\in\mathbb{N},n\le\omega_{N'}(x)\}\]
Let $P_{1}^{\epsilon}$ be defined as in the proof of theorem \ref{thm:stable}.
We can assume that $w_{p}(P_{1}^{\epsilon},\partial N')>2\epsilon$
for some for sufficiently small $\epsilon>0$. Because $A^{-}(N')=A(N')$
is convex, $\phi(\mu_{i})<\epsilon$ for $i=0,1$ implies \[
\phi(\mu_{t})<\epsilon.\]
Hence $P_{1}^{\epsilon}$ is convex and we can choose $B=P_{1}^{\delta}$
for some small $\delta>0$. Similarly $V(a)$ defined in theorem \ref{thm:cont}
is convex.

The set $N_{n}$ was defined as \[
N_{n}(\epsilon)=N\cap\operatorname{cl}\{y\,|\,\mbox{ \ensuremath{\exists x\in V(\epsilon)}, \ensuremath{m\ge0}\,\ s.t.\,\ \ensuremath{P_{n}^{[0,m]}(x)\subset\tilde{U}\,}and \ensuremath{P_{n}^{m}(x)=y}}\}\]
where $\tilde{U}=\operatorname{int}B$ and $N=\operatorname{cl}V(\epsilon_{0})$
are convex sets. Hence $N_{n}(\epsilon)$ is convex .\end{proof}
\begin{cor}
Under the assumption of the previous theorem if $p=1$ then all orbits
of $P_{n}$ are bounded for $n\ge n_{0}$, i.e. $P_{n}^{[0,\infty]}(\mu)\subset B_{R}^{w}(\mu_{0})$
for some $R\ge0$ and some fixed $\mu_{0}$. In particular, $K_{n}$
is the global weak attractor of $P_{n}$.\end{cor}
\begin{rem*}
The idea of the proof is to control the distance of $\mu_{1}$ and
$\mu_{0}$ by the distance of $\mu_{t}$ and $\mu_{1}$ where $\mu_{0}$
will be some stationary measure and $t\in(0,1]$ is sufficiently small.\end{rem*}
\begin{proof}
Using the Kantorovich-Rubinstein formula we have the following equality
for $\mu,\nu\in\mathcal{P}_{1}(X)$\[
w_{1}(\mu,\nu)=\inf_{\pi\in\Pi(\mu,\nu)}\int d(x,y)d\pi(x,y)=\sup_{\|\phi\|_{\operatorname{Lip}}\le1}\left\{ \int\phi d\mu-\int\phi d\nu\right\} ,\]
i.e. there is a sequence $\phi_{k}$ with $\|\phi_{k}\|_{\operatorname{Lip}}\le1$
such that $\int\phi_{k}d\mu-\int\phi_{k}d\nu\nearrow w_{1}(\mu,\nu)$.
Furthermore, there exist an optimal plan $\pi\in\Pi(\mu,\nu)$ such
that the infimum is actually attained.

Choose $\mu_{0}\in\mathcal{P}_{1}(X)$ and define $\mu_{t}=t\mu_{1}+(1-t)\mu_{0}$
for $t\in[0,1]$ and $\mu_{1}\in\mathcal{P}_{1}(X)$. We claim \[
w_{1}(\mu_{t},\mu_{0})=tw_{1}(\mu_{t},\mu_{0}).\]
Suppose $\pi\in\Pi(\mu_{1},\mu_{0})$ is the optimal plan and $\phi_{k}$
the sequence of Lipschitz maps as above. Then $\tilde{\pi}=t\pi+(1-t)(\operatorname{id},\operatorname{id})_{*}\mu_{0}$
is in $\Pi(\mu_{t},\mu_{0})$. Thus\[
w_{1}(\mu_{t},\mu_{0})\le\int d(x,y)d\tilde{\pi}(x,y)=tw_{1}(\mu_{1},\mu_{0}).\]
Furthermore, we have \[
w_{1}(\mu_{t},\mu_{0})\ge\int\phi_{k}d\mu_{t}-\int\phi_{k}d\mu_{0}=t\left(\int\phi_{k}d\mu_{1}-\int\phi_{k}d\mu_{0}\right).\]
Because the left hand side converges monotonically to $tw_{1}(\mu_{t},\mu_{0})$
we have proved our claim. 

Now fix some stationary measure $\mu_{0}\in K_{n}$. Since $N_{n}$
is a neighborhood of $K_{n}$ there is a $t\in(0,1]$ for all $\mu_{1}\in X$
such that $\mu_{t}$ as defined above is in $N_{n}$. Thus $P_{n}^{[0,\infty)}(\mu_{t})$
is in $N_{n}$ and bounded, i.e. $w_{1}(P_{n}^{k}(\mu_{t}),\mu_{0})\le R$
for some $R$. Because $P_{n}$ is convex linear and $\mu_{0}$ stationary
we have $P_{n}^{k}(\mu_{t})=tP_{n}^{k}(\mu_{1})+(1-t)\mu_{0}$ and
hence \begin{eqnarray*}
w_{1}(P_{n}^{k}(\mu_{1}),\mu_{0}) & = & \frac{1}{t}w_{1}(P_{n}^{k}(\mu_{t}),\mu_{0})\le\frac{R}{t},\end{eqnarray*}
which implies that the orbit of $\mu_{1}$ is bounded.
\end{proof}
\bibliographystyle{amsalpha}
\bibliography{ref}

\end{document}